\documentclass[12pt]{elsart}
\usepackage{amsfonts,amssymb}
\journal{Advances in Mathematics}

\newtheorem{Theorem}{Theorem}
\newtheorem{Definition}{Definition}
\newtheorem{Lemma}{Lemma}

\newenvironment{Proof}{{\bf Proof:}}{\hfill\rule{2mm}{2mm}\abst}

\newcommand{\abst}{\vspace{2ex}}

\begin{document}
\begin{frontmatter}
\title{Level Crossing Probabilities II: Polygonal Recurrence of 
Multidimensional Random Walks}
\author{Rainer Siegmund-Schultze and Heinrich von
Weizs\"{a}cker}
\address{Technische Universit\"{a}t Kaiserslautern\\
Fachbereich Mathematik\\
Erwin-Schr{\"o}dinger-Str., Geb\"{a}ude 48\\
67663 Kaiserslautern, Germany}

\begin{abstract}
In part I  we proved for an arbitrary one-dimensional random walk with
independent increments that the probability of crossing a level at a given
time $n$ is $O(n^{-1/2})$. In higher dimensions we call a random walk
'polygonally recurrent' if there is a bounded set, hit by infinitely many of the straight lines between two consecutive sites a.s. The above estimate implies that three-dimensional random
walks with independent components are polygonally transient. Similarly a
directionally reinforced random walk on $\mathbb{Z}^{3}$ in the sense of
Mauldin, Monticino and v.Weizs\"{a}cker \cite{Mauldin-Monticino-Weizsaecker}
is transient. On the other hand we construct an example of a transient but
polygonally recurrent random walk with independent components on $\mathbb{Z}%
^{2}$.
\end{abstract}

\end{frontmatter}
\maketitle

\section{Introduction}

This is a continuation of the paper \cite{SiegmundSchultze-Weizsaecker04a}
which gave a $O(n^{-1/2})$ bound for the level crossing probabilities of an
arbitrary one-dimensional random walk. We want to apply this result to study
'polygonal' transience and recurrence in higher dimensions and to
directionally reinforced random walks in the sense of \cite
{Mauldin-Monticino-Weizsaecker}.

\begin{Definition}
Let $(S_{n})$ be a random walk in $\mathbb{Z}^{d}$ or $\mathbb{R}^{d}$. We
call $(S_{n})$ \emph{polygonally recurrent} (resp. \emph{polygonally
transient}) if there is a bounded set $B$ (resp. there is no bounded set) such that a.s. there are infinitely many $n$ with the straight line between $S_{n},S_{n+1}$ hitting $B$.
\end{Definition}

A priori polygonal recurrence is a weaker statement than classical
recurrence, e.g. in one dimension every symmetric nontrivial random walk
oscillates between arbitrarily high negative and positive values and hence
is polygonally recurrent even if it is classically transient. In higher
dimensions it is less clear whether the two concepts really differ.

In three dimensions every (truly three-dimensional) random walk is transient.
If the components are independent then we get as a straightforward
consequence of our $O(n^{-1/2})$ estimate

\begin{Theorem}
\label{poltrans} A three-dimensional random walk whose three components are
independent is polygonally transient.
\end{Theorem}

We want to extend this result to the following situation, referred to as
''directionally reinforced random walk'' in \cite
{Mauldin-Monticino-Weizsaecker}: Let a particle move around in $\mathbb{Z}^d$
or $\mathbb{R}^d$. Assume that the particle moves with a constant velocity
along straight lines which are parallel to the coordinate axes, keeping its
direction of motion $\pm e_k,k\in \{1,...,d\}$ for a nonnegative finite
random time $T$ with $P(T>0)>0$ (in contrast to \cite{Mauldin-Monticino-Weizsaecker} we do not
require this random time to be strictly positive, cf. the first remark after
Theorem \ref{dirwalkex}), then changing to a different direction which is
chosen by equal chance among the $2d-1$ possibilities. This choice and the
random time spent to move in this new direction are assumed to be completely
independent of the past of the motion process. This process is, of course,
in general not a Markov process but, assuming that the first direction is
fixed, the successive locations of change into the first direction form a
truly $d-$dimensional random walk embedded in our process. In \cite
{Mauldin-Monticino-Weizsaecker} it was conjectured that in dimension 3 this
scheme is always transient in the sense that any bounded set is visited only
finitely often a.s. and that in dimension 2 the scheme is transient if the
embedded random walk is transient. (It is not difficult to see that for $d=1$
we have always recurrence and for $d>3$ always transience, \cite
{Mauldin-Monticino-Weizsaecker},  th. 3.1 and end of section 3). We prove in
section 2 the transience conjecture for $d=3$ using the above $O(n^{-1/2})$%
-bound.

However, we give in section 3 a somewhat involved example in 2 dimensions of
a directionally reinforced random walk which is recurrent whereas the
embedded random walk is transient but polygonally recurrent. Thus the level
crossing probabilities can be sufficiently higher than the return
probabilities to change a transience statement into recurrence.

\section{Transience in three dimensions}

Let us first give the simple\newline

\textbf{Proof of Theorem \ref{poltrans}.}  Let $A_n$ be the event that the
straight line between $S_n, S_{n+1}$ hits $[-1,1]^3$. Then by
our independence assumption 
\begin{eqnarray*}
P(A_n) \leq P(A_n^1)P(A_n^2)P(A_n^3)
\end{eqnarray*}
where $A_n^i$ denotes for $i = 1,2,3$ the event that the interval with the
endpoints $S^i_n, S^i_{n+1}$ meets $[-1,1]$, since $A_n$ implies each of the $A_n^i$. Clearly if $A_n^i$ occurs then
either $S_n^i \in [-1,1]$ or the the random walk $(S_n^i)$ crosses at time $n
$ at least one of the levels $-1$ or $1$. Both events have probability $O(n^{-1/2})$ by \cite{Spitzer}, p.72 (for the $\mathbb{Z}$--case) resp. \cite{Rosen}, Theorem 1 (for $\mathbb{R}$) and \cite
{SiegmundSchultze-Weizsaecker04a}, Theorem 2. Hence $P(A_n^i) = O(n^{-1/2})$
and $P(A_n) = O(n^{-3/2})$ which implies the result by Borel-Cantelli. \hfill%
\rule{2mm}{2mm}\vspace{2ex}

Consider now the model of $d$-dimensional directionally reinforced random
walk as it was defined in the introduction. We have, as a consequence of the
estimate for the probability of sign changes in symmetric one-dimensional
random walks, the transience of three-dimensional directionally reinforced
random walks. In \cite{Mauldin-Monticino-Weizsaecker} this was shown only
under the assumption that the waiting time between changes of direction has
a finite expectation and only for $d \ge 4$ without this moment condition.

\begin{Theorem}
\label{transdim3} For any dimension $d\geq 3$ the $d$-dimensional
directionally reinforced random walks are always transient in the sense that
bounded sets are visited only finitely often a.s.
\end{Theorem}

\begin{Proof}
1. Let us first modify the model to make the problem easier. Assume that,
when changing the direction of the travelling object, the next direction is
not chosen by equal chance from the $2d-1$ possible values which are
different from the previous one, but only from the $2d-2$ perpendicular
directions. We want to prove that a bounded set is visited only finitely
often. Fix a coordinate axis and call it 'vertical', and the others
'horizontal'. It is sufficient to show that the cube $[-1,1]^{d}$ is
penetrated or touched only finitely often, coming from vertical direction
(up or down). Consider those times when the particle changes from a
horizontal to vertical direction, or vice versa. Considering only these
times, the particle constantly changes from an independent symmetric random
walk in vertical direction $(\mathbf{R}^{1})$ to a horizontal symmetric and
independent random walk $(\mathbf{R}^{d-1})$ . Hitting the cube in the
assumed way means that the $\mathbf{R}^{d-1}-$random walk is just in the
cube $[-1,1]^{d-1}$, whereas the $\mathbf{R}^{1}-$random walk crosses the
levels $1$ or $-1$ or is in $[-1,1]$.

We have shown in \cite{SiegmundSchultze-Weizsaecker04a}, Theorem 2, that the probability of the second event is $O(n^{-\frac{1}{2}})$. The first event concerns a genuinely $d-1$-dimensional random walk and has a probability 
$O(n^{-(d-1)/2})$.
To see this we apply Theorem 3 of \cite{Esseen} which gives an estimate for the maximum probability of multi-dimensional rectangular domains with respect to the probability law of a sum of independent random vectors. In our case this estimate reads as
\begin{eqnarray*}
P(S_{n} &=&X_{1}+X_{2}+...X_{n}\in \lbrack -1,1]^{d-1}) \\
&\leq &C(\lambda )\left( 1-\sup_{x\in \mathbb{R}^{d-1}}P\left( X_{1}\in
D_{\lambda }+x\right) \right) ^{-\frac{d-1}{2}}n^{-\frac{d-1}{2}},
\end{eqnarray*}
for any $\lambda \leq 1$,
where $D_{\lambda }:=\{y=(y_{1},y_{2},...,y_{d-1})\in \mathbb{R}%
^{d-1}:|y_{j}|\leq \lambda $ for at least one $j\in \{1,2,...,d-1\}\}$. This
result is shown under a symmetry condition (S)\ meaning in our case that the
law of $X_{1}-X_{2}$ is invariant under any combined reflection of the
coordinate axes, i.e. under any application of a diagonal matrix of the form 
\begin{eqnarray*}
\left( 
\begin{array}{cccc}
\pm 1 & 0 & \cdots  & 0 \\ 
0 & \pm 1 &  & 0 \\ 
\vdots  &  & \ddots  & \vdots  \\ 
0 & 0 &  & \pm 1
\end{array}
\right).
\end{eqnarray*}
This is obviously fulfilled in the case considered here, since not only the
symmetrized law of $X_{1}-X_{2}$ but already the law of $X_{1}$ itself is
invariant under re-orientations of any axes. 

The only thing we have to prove is that $\sup_{x\in \mathbb{R}%
^{d-1}}P(X_{1}\in D_{\lambda }+x)\neq 1$ for a suitably chosen $\lambda $.
Due to construction, the law $P_{X_{1}}$of $X_{1}$ can be written as a
convex combination
$P_{X_{1}}=\delta Q_{1}+(1-\delta )Q_{2}$
with $Q_{1}$ being the law of a random vector $(\varepsilon
_{1}T_{1},\varepsilon _{2}T_{2},...,\varepsilon _{d-1}T_{d-1})$. Here the $%
\varepsilon _{i}$ and $T_{j}$ are completely independent of each other, the $%
T_{j}$ are i.i.d. distributed according to the time law of the directionally
reinforced random walk, and the $\varepsilon _{i}$ are i.i.d. coin tossing
random variables. This representation reflects the fact that with a positive
probability the moving object changes from the vertical motion to the first
horizontal direction, then to the second one and so on, and after that it
returns to vertical motion. Hence, it is sufficient to prove that $%
\sup_{x\in \mathbb{R}^{d-1}}Q_{1}(X_{1}\in D_{\lambda }+x)\neq 1$ for
sufficiently small $\lambda $. We have
\begin{eqnarray*}
\sup_{x\in \mathbb{R}^{d-1}}Q_{1}(X_{1} &\in &D_{\lambda }+x) \\
&=&\sup_{x\in 
\mathbb{R}^{d-1}}\left( 1-Q_{1}(|X_{1,j}-x_{j}|>\lambda
,j=1,2,...,d-1)\right)  \\
&=&\sup_{x\in \mathbb{R}^{d-1}}\left( 1-\prod_{j=1}^{d-1}\left(
1-Q_{1}(|X_{1,j}-x_{j}|\leq \lambda )\right) \right)  \\
&=&1-\left( 1-\sup_{t\in \mathbb{R}}P(|\varepsilon _{1}T_{1}-t|\leq \lambda
)\right) ^{d-1}.
\end{eqnarray*}
According to our assumption that $P(T>0)>0$, the law of $\varepsilon
_{1}T_{1}$ is non-degenerate. Hence the last expression is less than one for
sufficiently small $\lambda .$ Thus our first event has probability $O(n^{-(d-1)/2})$ and by Borel--Cantelli the modified model is transient for $d\geq 3$. 

2. Let us turn to the original model. The difference is that during a vertical 'phase' the
particle can change from up to down several times. Hence an intermediate
visit of $[-1,1]$ does not necessarily imply a level crossing, if we only
consider the positions at the beginning and the end of the vertical phase.

But, assuming infinitely many vertical visits to $[-1,1]^{d}$ with a
positive probability, we may consider the embedded process which, each time
the particle visits $[-1,1]^{d}$ during a vertical phase, registers
whether the following change takes it to a vertical direction or not. This
happens completely independently of the past with probabilities $p=\frac{1}{%
2d-1}$ and $q=\frac{2d-2}{2d-1}$, respectively. So the second case would
happen infinitely often, too, with the same positive probability. Hence also
this model would show infinitely many visits to $[-1,1]$ or crossings of levels $-1$ or $1$ of the
embedded vertical component during $[-1,1]^{d-1}-$visits of the horizontal
part. We may again apply Theorem 3 of \cite{Esseen} and Theorem 2 of \cite{SiegmundSchultze-Weizsaecker04a}
to disprove the possibility of such a behaviour.
\end{Proof}

\section{A two-dimensional example}

Despite the fact that return probabilities and level crossing probabilities
admit similar general asymptotic upper estimates nevertheless they can lead
to qualitatively different recurrence properties. We construct a transient
two-dimensional random walk which is 'polygonally recurrent' in a special
way.

\subsection{Results}

\begin{Theorem}
\label{Example} There is a symmetric distribution on the integers such that
two independent copies $(S_{n})$, $(\tilde{S}_{n})$ of the associated random
walk $(S_{n})$ satisfy the two conditions

(a) The two-dimensional random walk $(S_{n},\tilde{S}_{n})$ is transient.

(b) Almost surely the event $V_{n}=\{sgn(S_{n})=-sgn(S_{n+1}),\tilde{S}%
_{n}=0=\tilde{S}_{n+1}\}$ occurs for infinitely many $n$.

In particular

\begin{eqnarray}
&&\sum_{n=1}^{\infty }P(S_{n} =0)^{2}<\infty   \label{xconv}\\
&&\sum_{n=1}^{\infty }P(S_{n} =0)P(S_{n}S_{n+1}<0)=\infty .  \label{xdiv}
\end{eqnarray}
\end{Theorem}
Observe that due to (a) the event $V_{n}$ in (b) can be replaced by $\{S_{n}S_{n+1}<0,\tilde{S}%
_{n}=0=\tilde{S}_{n+1}\}$.
Moreover, note that in each term of the series (\ref{xdiv}) both factors are of the
order $O(n^{-1/2})$. Because of (\ref{xconv}) the first factor must be
actually of slightly smaller order, but the gap must be subtle because of
the divergence in (\ref{xdiv}). This lets us expect a somewhat delicate
construction. The construction will also yield a counterexample for a
(slightly modified, see the Remark below) conjecture on directionally
reinforced random walks in the sense of \cite{Mauldin-Monticino-Weizsaecker}:

\begin{Theorem}
\label{dirwalkex} There is a waiting time distribution on the nonnegative
integers such that the associated directionally reinforced random walk $%
(R_{m})$ on $\mathbb{Z}^{2}$ is recurrent but at any given lattice point the
walk a.s. changes direction only finitely often.
\end{Theorem}

\noindent \textbf{Remark} The actual waiting time distribution constructed
below gives positive probability to the value $0$. One could insist on a
strictly positive waiting time in order to be exactly in the framework of 
\cite{Mauldin-Monticino-Weizsaecker} but the conceptual arguments based on
unimodality below would not be directly applicable.

\subsection{The main idea of the proof}

Intuitively the construction of Theorem \ref{Example} is based on the
observation that if the one-dimensional random walk $(S_{n})$ has lattice
constant $1$ and the underlying symmetric random variable $X$ has finite
variance $\sigma ^{2}$ then 
\begin{equation}
P(S_{n}=0)\asymp \frac{1}{\sigma \sqrt{n}}  \label{vconv}
\end{equation}
and if $\gamma _{\alpha }$ is a $(1-\alpha )$-quantile of the distribution
of $|X|$ then 
\begin{equation}
P(S_{n}S_{n+1}<0)\geq \frac{\alpha }{2}P(|S_{n}|\leq \gamma _{\alpha
})\asymp \frac{\alpha \gamma _{\alpha }}{\sigma \sqrt{n}}.  \label{vdiv}
\end{equation}
So for (\ref{xconv}) the variance must be infinite. But we can let the high
values of $X$ occur rarely enough so that with very high probability the
behaviour of our random walk up to time $n$ equals the behaviour of another
random walk with variance $\sigma _{n}^{2}$ and corresponding quantiles $%
\gamma _{\alpha _{n}}$ of the absolute value where these numbers grow in
such a balanced manner that if we plug them into (\ref{vconv}) and (\ref
{vdiv}) we get (\ref{xconv}) and (\ref{xdiv}). Clearly (\ref{xconv}) then
implies the transience (a). If the events $V_{n}$ would be independent, one would easily infer statement (b) from (\ref{xdiv}). Since they are not, an extra argument is needed.
The key to this step is Lemma \ref{recurrevents} (see Appendix) which gives a
quantitative version of the fact known eg from Markov chains that the
frequency of certain events is high \textit{with high probability} if only
its \textit{expectation} is high enough.

For the estimates in the main body of the proof it is useful to have
symmetric unimodal distributions (since the notion of unimodality is used in the literature not completely consistently, see the Appendix for a definition). 

The waiting time distribution in Theorem \ref{dirwalkex} will be given as a
mixture of uniform distributions 
\begin{equation}
\mathcal{L}(T)=\sum_{l=1}^{\infty }p_{l}R[0,y_{l}],\;\sum_{l=1}^{\infty
}p_{l}=1\;  \label{Tdef}
\end{equation}
where $R[0,y_{l}]$ denotes the uniform distribution on an integer interval $%
[0,y_{l}]$ and the increasing resp. decreasing sequences $(y_{l})$ and $%
(p_{l})$ will be constructed recursively below. Observe that $T$ is a random
variable with non-increasing weights, i.e. $P(T=k)\geq P(T=k+1)$ for each $%
k\in \{0,1,2,...\}.$

In the following we will make use of a coupled sequence of random walks with finite
variance waiting times which approximates our infinite variance random walk
given by the above waiting time distribution. For this end we first define a
one-sided sequence $\mathbf{T}=\{T,\kappa ,T^{(1)},T^{(2)},...\}$ of random
variables with the property $T^{(\kappa )}=T^{(\kappa +1)}=...=T$ for the
random index $\kappa $ and such that for all $m > k$ we have 
\begin{equation}  \label{tkindep}
\mathcal{L}(T^{(k)}|\kappa = m) = \mathcal{L}(T^{(k)})=z_{k}^{-1}%
\sum_{l=1}^{k}p_{l}R[0,y_{l}],\ z_{k}=\sum_{l=1}^{k}p_{l}.
\end{equation}
To do this, think of $T$ as the result of a two-step construction: First the
probability distribution $\{p_{l}\}$ is used to find a random index $\kappa $%
, then $T$ is realized by choosing a random integer from $[0,y_{\kappa }]$
according to the uniform distribution on this interval. Given $(T,\kappa ),$
choose $\{T^{(1)},T^{(2)},...,T^{(\kappa -1)}\}$ as a sequence of
independent realizations of $\mathcal{L}(T^{(i)}),i<\kappa ,$ respectively,
and let $T^{(\kappa )}=T^{(\kappa +1)}=...=T.$

This leads to a hierarchical structure: The random walk will be a certain
mixture of a sequence of random walks on different scales of time and space.
The walk at level $k$ runs on a space scale determined by $y_{k}$ and a step
frequency determined by $p_{k}$. The reader should think of really rapidly
increasing scales. In fact the simplest lower estimates for which the
construction below (cf. (\ref{ck})-(\ref{lest1})) can be carried out show
that the $y_{k}$ grow at least like in a recursion of the form 
\begin{eqnarray*}
y_{k+1}=e^{const\cdot y_k}.
\end{eqnarray*}

The main part of the proof of Theorem \ref{Example} lies in the careful choice of the
parameters $p_k$ and $y_k$ of the waiting time $T$ given by (\ref{Tdef}). Once they are determined we choose an iid sequence $(T_i)$ with the same law as $T$ and consider\newline

i) The directionally reinforced random walk $(R_m)$ on $\mathbb{Z}^2$ which
moves in unit size steps and starts at the origin horizontally in either the
positive or the negative direction. After the waiting time $T_1$ it switches
with uniform probability to one of the three other directions. After waiting
time $T_2$ it changes direction again and so on. We want to show that for
our particular law of waiting times a.s. $R_m$ visits the origin infinitely
often but it changes direction at the origin only finitely many times. This
will prove Theorem \ref{dirwalkex}.\newline

ii) The sequence $(S_n,\tilde{S}_n)$ consisting of the successive locations
at which $(R_m)$ changes from vertical movement back to horizontal movement.
By the properties of $R_m$ the increments of $(S_n,\tilde{S}_n)$ are iid.
with independent components. In fact the law of $S_n - S_{n-1}$ is equal to
the law of the random variable $X$ defined in Lemma \ref{unimod} where $%
\epsilon$ determines the sign of the first part of the horizontal movement
of $R_m$ after the visit of $(S_{n-1},\tilde{S}_{n-1})$ and $G$ is a
geometric random variable with parameter $2/3$ which determines the number
of horizontal flips before the next vertical step at the location $(S_n,%
\tilde{S}_{n-1})$. Similarly the second component $\tilde{S}_n - \tilde{S}%
_{n-1}$ of the increment has the same law and it determines the following
vertical movement from $(S_n,\tilde{S}_{n-1})$ to $(S_n,\tilde{S}_n)$.

Thus by Lemma \ref{unimod} $(S_n,\tilde{S}_n)$ is a random walk with
independent components which have a symmetric law as required in Theorem \ref
{Example}. Now it is not hard to see that the assertion of Theorem \ref{Example} implies Theorem \ref{dirwalkex}. In fact, consider first the conditional probability $z_{r,t}$ that the directionally reinforced random walk considered above never again after time $t$ changes direction at the origin, given that at time $t$ it changes direction at the origin, coming from direction $r$. Obviously by construction this probabilty does not depend on $r,t$, so we denote it by $z$. At each instance where the walk changes direction at the origin, with probability $z$ it will never do so again, completely independent from the past. So in the case $z>0$ there will be a.s. only finitely many changes of direction at the origin, while for $z=0$ there will be a.s. infinitely many such events. Now assume the latter, i.e. $z=0$. Note that at each change of direction at the origin the following event is independent of the past and has positive probability: The walk changes to a perpendicular direction, makes zero steps in this direction, then changes again to a perpendicular direction. By assumption $z=0$ this will happen a.s. infinitely often, too. But by construction of the embedded random walk $(S_n,\tilde{S}_n)$ this means that $(S_n,\tilde{S}_n)$ visits the origin infinitely often a.s. Hence the transience part (a) of Theorem \ref{Example} implies $z>0$ which means that there are a.s. only finitely many changes of direction at the origin for the directionally reinforced random walk. The same argument applies to any other lattice point.
On the other hand, part (b) of Theorem \ref{Example} immediately implies that there are a.s. infinitely many visits of the origin for the directionally reinforced random walk $(R_m)$.
Thus it suffices to find a waiting
time distribution ensuring the validity of Theorem \ref{Example}.

\subsection{The hierarchical construction}

In this section the underlying parameters $p_{k}$ and $y_{k}$ of the the law
of the waiting time (\ref{Tdef}) are not yet fixed.

We start with an i.i.d. two-dimensional array 
\begin{eqnarray*}
(\mathbf{T}_{i,n})=(\{T_{i,n},\kappa
_{i,n},T_{i,n}^{(1)},T_{i,n}^{(2)},...\})
\end{eqnarray*}
of random sequences as constructed above. Let $(\epsilon _{n})$ and $(G_{n})$
be two iid sequences of cointossing resp. geometric (with parameter $2/3$)
random variables chosen independently from $(\mathbf{T}_{i,n})$ and let 
\begin{eqnarray}
X_{n} &=&\epsilon _{n}\sum_{i=1}^{G_{n}}(-1)^{i}T_{i,n},  \label{Xkdef} \\
X_{n}^{(k)} &=&\epsilon _{n}\sum_{i=1}^{G_{n}}(-1)^{i}T_{i,n}^{(k)}.
\end{eqnarray}
Then we can define the two-dimensional random walks $(S_{n},\tilde{S}%
_{n})$ and $(S_{n}^{(k)},\tilde{S}%
_{n}^{(k)})$ by 
\begin{eqnarray}
S_{n} &=&\sum_{j=1}^{n}X_{j},\ \ \tilde{S}_{n}=\sum_{j=1}^{n}%
\widetilde{X}_{j},  \nonumber \\
S_{n}^{(k)} &=&\sum_{j=1}^{n}X_{j}^{(k)},\ \ \tilde{S}_{n}^{(k)}=%
\sum_{j=1}^{n}\widetilde{X}_{j}^{(k)},
\end{eqnarray}
where we take $(\widetilde{\mathbf{T}}_{i,n})=(\{\widetilde{T}_{i,n},%
\widetilde{\kappa }_{i,n},\widetilde{T}_{i,n}^{(1)},\widetilde{T}%
_{i,n}^{(2)},...\}),(\widetilde{\varepsilon }_{n})$, and $(\widetilde{G}_{n})
$ as above independent of $(\mathbf{T}_{i,n}),(\varepsilon _{n})$, and $%
(G_{n})$, and then use them to construct $\widetilde{X}_{n},\widetilde{X}%
_{n}^{(k)},\widetilde{S}_{n}$ and $\widetilde{S}_{n}^{(k)}$ as we
constructed $X_{n},X_{n}^{(k)},S_{n}$ and $S_{n}^{(k)}$. 

In the sequel, a truncated version of this coupled sequence of
random walks will be useful. It is obtained by cutting the random variable $%
\kappa $ defining the hierarchy level at a given value $K$. Hence we define
a truncated version of the sequence $\mathbf{T}$ by $\mathcal{L}(\mathbf{T}%
^{|K})=\mathcal{L}(\mathbf{T}|\kappa \leq K)$, i.e. $\mathbf{T}%
^{|K}=\{T^{|K},\kappa ^{|K},T^{(1)|K},T^{(2)|K},...,T^{(K)|K}\}$ with $%
\mathcal{L}(T^{|K})=\mathcal{L}(T^{(K)}),\mathcal{L}(\kappa ^{|K})=\mathcal{L%
}(\kappa |\kappa \leq K)=\{p_{l}/z_{K}\}_{l=1}^{K}$ and $\mathcal{L}%
(T^{(k)|K})=\mathcal{L}(T^{(k)})$. Observe that this truncation shares with
the original coupling the property that $T^{(1)|K},T^{(2)|K},...,T^{(\kappa
^{|K})|K}$ are (conditionally with respect to $\kappa ^{|K}$) independent,
while $T^{(\kappa ^{|K})|K},T^{(\kappa ^{|K}+1)|K},...,T^{(K)|K}$ coincide.

So the truncated version of the coupled sequence
of random walks is obtained by substituting the i.i.d. sequence $(\mathbf{T}_{i,n})$
with the sequence $(\mathbf{T}_{i,n}^{|K})$. We define
\begin{eqnarray}  
X_{n}^{(k)|K} &=&\epsilon _{n}\sum_{i=1}^{G_{n}}(-1)^{i}T_{i,n}^{(k)|K}. \label{truncb} 
\end{eqnarray}
This yields a
finite sequence of random walks $\{S_{n}^{(k)|K},\widetilde{S}%
_{n}^{(k)|K}\}_{k=1}^{K}$ with the property 
\begin{equation}
\mathcal{L}(S_{n}^{(k)|K})=\mathcal{L}(S_{n}^{(k)}),k\leq K.  \label{trunca}
\end{equation}
 
For the transience proof we prepare now an upper
estimate of the return probabilities of $S_{n}^{(k)}$.
For this end we use the following estimate

\begin{Lemma}
There is a positive constant $A$ such that the return probabilities $%
s_{n}^{k}=P(S_{n}^{(k)}=0)$ satisfy the recursive estimate 
\begin{equation}
s_{n}^{k}\leq s_{n}^{k-1}(1-p_{k})^{n}+\frac{A}{p_{k}y_{k}}\frac{1}{\sqrt{n}}%
.  \label{qkn2}
\end{equation}
\end{Lemma}

\begin{Proof}

We set the truncation parameter $K=k$ and denote by $Z_{n}^{(k)}$ the number of pairs $(i,j)$ with $j\leq n$ and $%
i\leq G_{j}$ such that in the representation (\ref{truncb}) the r.v. $T^{(k)|k}_{i,j}$
actually is of the maximal level $k$, i.e. $\kappa^{|k} _{i,j}= k.$ For each
index pair this happens with probability $p_{k}/z_{k}>p_{k}$.

On the set $\{Z_{n}^{(k)}=0\}$ we have $S_{n}^{(k)|k}=S_{n}^{(k-1)|k}$ by
construction, hence by (\ref{trunca})
\begin{eqnarray*}
P(S_{n}^{(k)} =0)&=&P(S_{n}^{(k)|k}=0) \\
&=&P(S_{n}^{(k-1)|k}=0)P(Z_{n}^{(k)}=0|S_{n}^{(k-1)|k}=0) \\
& &+P(S_{n}^{(k)|k}=0,Z_{n}^{(k)}>0)
\\
&\leq &P(S_{n}^{(k-1)}=0)P(\kappa _{1,j}^{|k}<k\textrm{ for }j\leq
n|S_{n}^{(k-1)|k}=0)\\
& &+P(S_{n}^{(k)|k}=0,Z_{n}^{(k)}>0).
\end{eqnarray*}
Our whole construction is designed to ensure that the conditional law of $\{S_{n}^{(k-1)|k}=0\}$ given $\{\kappa _{i,j}^{|k}=k\}$ is equal to the unconditional law for any $i$ and $j\leq n$. Therefore, the two events $\{\kappa _{1,j}^{|k}<k$ for $j\leq n\}
$ and $\{S_{n}^{(k-1)|k}=0\}$ are independent of each other, and the
probability of the first event is $(1-p_{k}/z_{k})^{n}\leq (1-p_{k})^{n}$.
Consequently, 
\begin{eqnarray}
P(S_{n}^{(k)} =0)&\leq&
P(S_{n}^{(k-1)|k}=0)(1-p_{k})^{n}+P(S_{n}^{(k)|k}=0,Z_{n}^{(k)}>0) \label{uestkzerl} \\
&=&P(S_{n}^{(k-1)}=0)(1-p_{k})^{n}+P(S_{n}^{(k)|k}=0,Z_{n}^{(k)}>0). \nonumber 
\end{eqnarray}
In order to estimate the second term we note that conditionally on the
knowledge of the set of pairs $(i,j)$ with a contribution of level $%
k,$ and of the signs $\epsilon _{j}$
the law of the sum $S_{n}^{(k)|k}$ is of the form as in (\ref{maxest}) with $%
y = y_{k}$, except for an additional convolution factor (coming from the contribution of terms of level less than $k$) which does not increase the maximum probability. An upper estimate of this maximum probability is conserved under convex
combinations and hence the maximum probability of the conditional law of $S_{n}^{(k)|k}$
given the value of $Z_{n}^{(k)}$ is at most $\frac{D}{y_{k}\sqrt{Z_{n}^{(k)}}%
}$. Thus 
\begin{eqnarray}
P(S_{n}^{(k)|k}=0,Z_{n}^{(k)}>0) &=&\sum_{m=1}^{\infty
}P(S_{n}^{(k)|k}=0|Z_{n}^{(k)}=m)P(Z_{n}^{(k)}=m)  \label{uestcond} \\
&\leq &\sum_{m=1}^{\infty }\frac{D}{\sqrt{m}y_{k}}P(Z_{n}^{(k)}=m)  \nonumber \\
&=&\frac{D}{y_{k}}\mathbb{E}\left( \frac{1}{\sqrt{Z_{n}^{(k)}}}\mathbf{1}%
_{\{Z_{n}^{(k)}>0\}}\right) .  \nonumber
\end{eqnarray}
Clearly $\mathbb{E}(Z_{n}^{(k)})=np_{k}/z_{k}\mathbb{E}(G)$ and by Wald's
identity \textrm{Var}$(Z_{n}^{(k)})=np_{k}/z_{k}(1-p_{k}/z_{k})\mathbb{E}%
(G)+n(p_{k}/z_{k})^{2}\textrm{Var}(G)$. Thus Chebyshev's inequality gives two positive constants $b,B$ such
that 
\begin{eqnarray}
P(bnp_{k}/z_{k}\leq Z_{n}^{(k)})&\geq& 1-\frac{np_{k}/z_{k}(1-p_{k}/z_{k})\mathbb{E}(G)+n(p_{k}/z_{k})^{2}\textrm{Var}(G)}{((\mathbb{E}(G)-b)np_{k}/z_{k})^{2}} \nonumber \\
&\geq& 1-\frac{B}{%
np_{k}}.  \label{ChZ}
\end{eqnarray}
Therefore for some positive constant $A$ 
\begin{equation}
\mathbb{E}\left( \frac{1}{\sqrt{Z_{n}^{(k)}}}\mathbf{1}_{\{Z_{n}^{(k)}>0\}}%
\right) \leq \frac{1}{\sqrt{bnp_{k}/z_{k}}}+\frac{B}{np_{k}}\leq 
\frac{A}{p_{k}\sqrt{n}}.  \label{Ewurz}
\end{equation}
Plugging the estimates (\ref{uestcond}) and (\ref{Ewurz}) into (\ref
{uestkzerl}) yields the desired result.
\end{Proof}

\subsection{Recursive choice of the parameters}

We start with $y_{1}=1$ and $p_{2}=\frac{1}{4}$. The quantity $p_{1}$ will
be chosen only in the end in order to get a total sum $1$, but with
condition (\ref{pkhalf}) below it is obvious that $\sum_{k=2}^{\infty
}p_{k}\leq \frac{1}{2}$ and hence $p_{1}$ is at least 1/2. Let now $k\in 
\mathbb{N},k>1$ be given and assume that $y_{l}$ for $1\leq l<k$ and $p_{l}$
for $2\leq l\leq k$ are already defined. Then choose an integer $c_{k}$ with 
\begin{equation}
c_{k}>\frac{k^{8}}{p_{k}^{2}}.  \label{ck}
\end{equation}
Now choose the two numbers $y_{k},p_{k+1}$ such that 
\begin{eqnarray}
y_{k}\sqrt{p_{k}} &\geq &\max (12c_{k},y_{k-1}\sqrt{p_{k-1}})  \label{pest}
\\
0 &< & p_{k+1} \leq \frac{1}{2}p_{k}  \label{pkhalf}
\end{eqnarray}
and 
\begin{equation}
\frac{1}{2k^{4}}\leq (\frac{A}{p_{k}y_{k}})^{2}\log \frac{1}{p_{k+1}}\leq 
\frac{1}{k^{4}}. \label{lest1}
\end{equation}
Observe that we may guarantee (\ref{pest}) to hold for $k=2$, even though $p_{1}$ is unknown in the beginning.
This completes the recursive construction.

\subsection{Transience of $(S_n,\tilde{S}_n)$}

We now have to verify that the resulting distribution has the desired
properties. Let $(\tilde{S}_{n})$ be an independent copy of the random walk $%
(S_{n})$. 
For the transience we prove the convergence of the series (\ref{xconv}). 
By construction, the indicator function of the event $S_{n}=0$ is the pointwise limit of the indicator function of $S_{n}^{(k)}=0$.
Hence, for a fixed integer $N$ we get 
\begin{eqnarray*}
\sum_{n=1}^{N}P(S_{n}=0,\tilde{S}_{n}=0) &=
&\sum_{n=1}^{N}\lim_{k\rightarrow \infty }P(S_{n}^{(k)}=0,\tilde{S}%
_{n}^{(k)}=0)(1-p_{k+1})^{2n} \\
&\leq &\lim_{k\rightarrow \infty }\sum_{n=1}^{\infty
}(s_{n}^{k})^{2}(1-p_{k+1})^{2n}
\end{eqnarray*}
The proof of (a) will be complete if we can verify that for each $k$ 
\begin{equation}
\sqrt{\sum_{n=1}^{\infty }(s_{n}^{k})^{2}(1-p_{k+1})^{2n}}\leq
2\sum_{j=1}^{k}\frac{1}{j^{2}}  \label{pi}
\end{equation}
because this will show that the series in (\ref{xconv}) is $\leq \frac{\pi ^{4}}{%
9}$.

We prove (\ref{pi}) by induction. In the case $k=1$ this follows from $%
s_{n}^{1}\leq 1$ and $p_{2}=\frac{1}{4}$. For the induction step observe 
\begin{equation}
\sum_{n=1}^{\infty }\frac{1}{n}(1-q)^{n}=-\log (1-(1-q))=\log (1/q)
\label{logreihe}
\end{equation}
for $0<q<1$. We use (\ref{qkn2}) and the triangle inequality in the sequence
space $\ell ^{2}$ and get 
\begin{eqnarray*}
&&\sqrt{\sum_{n=1}^{\infty }(s_{n}^{k})^{2}(1-p_{k+1})^{2n}} \\
&\leq &\sqrt{\sum_{n=1}^{\infty
}(s_{n}^{k-1})^{2}(1-p_{k})^{2n}(1-p_{k+1})^{2n}}+\sqrt{(\frac{A}{p_{k}y_{k}}%
)^{2}\sum_{n=1}^{\infty }\frac{1}{n}(1-p_{k+1})^{n}} \\
&\leq &\sqrt{\sum_{n=1}^{\infty }(s_{n}^{k-1})^{2}(1-p_{k})^{2n}}+\sqrt{(%
\frac{A}{p_{k}y_{k}})^{2}\log \frac{1}{p_{k+1}}} \\
&\leq &2\sum_{j=1}^{k-1}\frac{1}{j^{2}}+\sqrt{\frac{1}{k^{4}}}.
\end{eqnarray*}
In the last step we have used the induction hypothesis for the first term
and (\ref{lest1}) for the second term.

This proves (\ref{pi}) and hence part (a) of Theorem \ref{Example}. 

\subsection{Recurrence of the events $V_n$}
Of course we want to compare our two-dimensional random walk with
its approximations. Therefore the sets 
\begin{eqnarray*}
F_{n,k}=\{\kappa _{i,j}\leq k,1\leq j\leq n,1\leq i\leq G_{j}\}\cap \{\widetilde{%
\kappa }_{i,j}\leq k,1\leq j\leq n,1\leq i\leq \widetilde{G}_{j}\}
\end{eqnarray*}
are important.
We consider the events 
\begin{eqnarray*}
E_{n,k}=\{\tilde{S}_{n}=0,|S_{n}|\leq c_{k}\}\cap F_{n,k}.
\end{eqnarray*}
We introduce the notation $p_{k}^{\ast }:=\sum_{l \geq k}p_{k} \geq p_{k}$.
We write $\mathbf{G}$ for the $n$-tuple $(G_{1},G_{2},...,G_{n})$ and $%
\widetilde{\mathbf{G}}$ for the $n$-tuple $(\widetilde{G}_{1},\widetilde{G}%
_{2},...,\widetilde{G}_{n})$. By construction, conditioned on $\mathbf{G}$
the random variable $S_{n}^{(k)}$ is independent of $F_{n,k}$ and on the set 
$F_{n,k}$ we have $S_{n}^{(k)}=S_{n}$. \ Hence
\begin{equation}
P(E_{n,k})=\mathbb{E}(P(F_{n,k}|\mathbf{G,}\widetilde{\mathbf{G}}%
)P(|S_{n}^{(k)}|\leq c_{k}|\mathbf{G})P(\widetilde{S}_{n}^{(k)}=0|\widetilde{%
\mathbf{G}}))
\label{asymkr}.
\end{equation}
Moreover, conditioning on $\mathbf{G}$, Lemma \ref{unimod} gives that the laws of $%
X_{1}^{(k)},...,X_{n}^{(k)}$ are symmetric unimodal and so is the
conditional law of $S_{n}^{(k)}$. Also by Lemma \ref{unimod} we have
\begin{eqnarray*}
\textrm{Var}(S_{n}^{(k)}|\mathbf{G})=\sum_{1\leq i\leq n}\textrm{Var}%
(X_{i}^{(k)}|G_{i})\leq 4\sum_{1\leq i\leq n}G_{i}\textrm{Var}(T_{1,1}^{(k)})
\end{eqnarray*}
and
\begin{eqnarray*}
\textrm{Var}(S_{n}^{(k)}|\mathbf{G})=\sum_{1\leq i\leq n}\textrm{Var}%
(X_{i}^{(k)}|G_{i})\geq n\textrm{Var}(T_{1,1}^{(k)}).
\end{eqnarray*}
Corresponding relations are valid for $\widetilde{S}_{n}^{(k)}$.
From (\ref{tkindep}) we get $\mathbb{E}((T_{1}^{(k)})^{2})\geq \frac{%
1}{4}p_{k}y_{k}^{2}$ and hence by Lemma \ref{momest} $\mathrm{Var}%
(T_{1,1}^{(k)})\geq \frac{1}{12}p_{k}y_{k}^{2}$ for large enough $k$.
Hence Var$(S_{n}^{(k)}|\mathbf{G})\geq np_{k}y_{k}^{2}/12$, and this
expression is at least $12c_{k}^{2}$ $\geq 12$ by (\ref{pest}). Now Lemma \ref{unimodest} can
be applied to get
\begin{eqnarray*}
P(\widetilde{S}_{n}^{(k)} =0|\widetilde{\mathbf{G}})&\geq& \frac{d}{\sqrt{%
\textrm{Var}(S_{n}^{(k)}|\mathbf{G})}} \\
P(|S_{n}^{(k)}| \leq c_{k}|\mathbf{G})&\geq& \frac{dc_{k}}{\sqrt{\textrm{Var}%
(S_{n}^{(k)}|\mathbf{G})}}.
\end{eqnarray*}
Hence
\begin{eqnarray*}
P(|S_{n}^{(k)}|\leq c_{k}|\mathbf{G})P(\widetilde{S}_{n}^{(k)}=0|\widetilde{%
\mathbf{G}})\geq \frac{d^{2}c_{k}}{4\textrm{Var}(T_{1,1}^{(k)})\sqrt{\left(
\sum_{1\leq i\leq n}G_{i}\right) \left( \sum_{1\leq i\leq n}\widetilde{G}%
_{i}\right) }}.
\end{eqnarray*}
We have $P(F_{n,k}|\mathbf{G,}\widetilde{\mathbf{G}})=(1-p_{k+1}^{\ast
})^{\sum_{1\leq i\leq n}G_{i}+\widetilde{G}_{i}}$, hence from (\ref{asymkr}) we get
\begin{eqnarray*}
P(E_{n,k})\geq \frac{d^{2}c_{k}}{4\textrm{Var}(T_{1,1}^{(k)})}\left( \mathbb{E}%
\frac{(1-p_{k+1}^{\ast })^{\sum_{1\leq i\leq n}G_{i}}}{\sqrt{\sum_{1\leq
i\leq n}G_{i}}}\right) ^{2}.
\end{eqnarray*}
We consider the function $\psi _{a}(\lambda ):=\exp (-a\lambda )/\sqrt{%
\lambda },a>0$. It is easily checked that this function is convex for $%
\lambda >0$. Hence by Jensen's inequality we get
\begin{eqnarray*}
P(E_{n,k})\geq \frac{d^{2}c_{k}}{4\textrm{Var}(T_{1,1}^{(k)})}\frac{1}{n%
\mathbb{E}G}(1-p_{k+1}^{\ast })^{2n\mathbb{E}G}.
\end{eqnarray*}
We have the estimate 
\begin{eqnarray*}
\mathrm{Var}(T_{1,1}^{(k)})\leq \mathbb{E}((T_{1,1}^{(k)})^{2})\leq
z_{k}^{-1}\sum_{l=1}^{k}p_{l}y_{l}^{2}\leq z_{k}^{-1}kp_{k}y_{k}^{2}
\end{eqnarray*}
where the last inequality follows from (\ref{pest}).
So we arrive at 
\begin{eqnarray*}
P(E_{n,k})\geq \frac{d^{2}c_{k}}{4z_{k}^{-1}kp_{k}y_{k}^{2}}\frac{1}{n%
\mathbb{E}G}(1-p_{k+1}^{\ast })^{2n\mathbb{E}G}.
\end{eqnarray*}
Consequently, since $p_{k+1}^{\ast }\leq 2p_{k+1}$ by (\ref{pkhalf})  
\begin{eqnarray*}
\sum_{n=1}^{\infty }P(E_{n,k}) &\geq &\frac{d^{2}c_{k}}{\mathbb{E}%
G4z_{k}^{-1}kp_{k}y_{k}^{2}}\sum_{n=1}^{\infty }\frac{1}{n}(1-p_{k+1}^{\ast
})^{2n\mathbb{E}G} \\
&\geq &\frac{d^{2}c_{k}}{\mathbb{E}G4z_{k}^{-1}kp_{k}y_{k}^{2}}%
\sum_{n=1}^{\infty }\frac{1}{n}(1-2p_{k+1})^{2n\mathbb{E}G} \\
&\geq &\frac{d^{2}c_{k}}{\mathbb{E}G4z_{k}^{-1}kp_{k}y_{k}^{2}}%
\sum_{n=1}^{\infty }\frac{1}{n}(1-4\mathbb{E}Gp_{k+1})^{n} \\
&=&\frac{d^{2}c_{k}}{\mathbb{E}G4z_{k}^{-1}kp_{k}y_{k}^{2}}\log \frac{1}{4%
\mathbb{E}Gp_{k+1}}.
\end{eqnarray*}

If $\Phi _{k}$ denotes the total number of the events $E_{n,k}$ which occur
we get from the conditions (\ref{ck}) and (\ref{lest1}) for large enough $k$ 
\begin{equation}
\mathbb{E}(\Phi _{k})>\frac{k^{2}}{p_{k}}.  \label{expcount}
\end{equation}
Next we want to apply Lemma \ref{recurrevents}. Let $\mathcal{F}_{m}$ denote the $\sigma $-field generated by the random variables $%
S_{j},S_{j}^{(k)},G_{j},\kappa _{i,j}$ and $\tilde{S}_{j},\tilde{S}%
_{j}^{(k)},\tilde{G}_{j},\widetilde{\kappa }_{i,j}$ with $1\leq j\leq m,1\leq i.$ 
Set $\mathbf{G}_{l}^{m}=(G_{l},G_{l+1},...,G_{m}),\widetilde{\mathbf{G}}%
_{l}^{m}=(\widetilde{G}_{l},\widetilde{G}_{l+1},...,\widetilde{G}_{m})$ By
conditional independence on $E_{m,k}$ we have the relation
\begin{eqnarray*}
P(E_{n,k}|\mathcal{F}_{m}) &=&
\mathbb{E}(P(F_{n,k}|F_{m,k},\mathbf{G}_{m+1}^{n},\widetilde{\mathbf{G}}%
_{m+1}^{n}) \\
& & \cdot P(|S_{n}^{(k)}| \leq c_{k}|S_{m}^{(k)},\mathbf{G}_{m+1}^{n})P(|%
\widetilde{S}_{n}^{(k)}|=0|\widetilde{S}_{m}^{(k)},\widetilde{\mathbf{G}}%
_{m+1}^{n})).
\end{eqnarray*}
For any $n-m$-tupel of positive integers $\mathbf{a}%
=(a_{1},a_{2},...,a_{n-m})$ and any integer $x$ we obtain
\begin{eqnarray*}
P(|S_{n}^{(k)}| \leq c_{k}|S_{m}^{(k)}=x,\mathbf{G}_{m+1}^{n}=\mathbf{a}%
)&=&P(|S_{n-m}^{(k)}-x|\leq c_{k}|\mathbf{G}_{1}^{n-m}=\mathbf{a}) \\
&\leq &P(|S_{n-m}^{(k)}|\leq c_{k}|\mathbf{G}_{1}^{n-m}=\mathbf{a})
\end{eqnarray*}
since the conditional law $P(S_{n-m}^{(k)}\in (\cdot ) |\mathbf{G}_{1}^{n-m}=%
\mathbf{a})$ is unimodal symmetric, and for the same reason we get
\begin{eqnarray*}
P(\widetilde{S}_{n}^{(k)}=0|\widetilde{S}_{m}^{(k)}=x,\widetilde{\mathbf{G}}%
_{m+1}^{n}=\mathbf{a})\leq P(\widetilde{S}_{n-m}^{(k)}=0|\widetilde{\mathbf{G%
}}_{1}^{n-m}=\mathbf{a}).
\end{eqnarray*}
Let $\mathbf{b}=(b_{1},b_{2},...,b_{n-m})$ another $n-m$-tupel of positive
integers. We get on $E_{m,k}$%
\begin{eqnarray*}
P(F_{n,k}|F_{m,k},\mathbf{G}_{m+1}^{n}=\mathbf{a},\widetilde{\mathbf{G}}%
_{m+1}^{n}=\mathbf{b})=P(F_{n-m,k}|\mathbf{G}_{1}^{n-m}=\mathbf{a},%
\widetilde{\mathbf{G}}_{1}^{n-m}=\mathbf{b}).
\end{eqnarray*}
So %
\begin{eqnarray}
P(E_{n,k}|\mathcal{F}_{m}) &\leq&  \nonumber \\ 
\mathbb{E}(P(F_{n-m,k}&|&\mathbf{G}_{1}^{n-m},\widetilde{\mathbf{G}}%
_{1}^{n-m})P(|S_{n-m}^{(k)}| \leq c_{k}|\mathbf{G}_{1}^{n-m})P(|\widetilde{%
S}_{n-m}^{(k)}|=0|\widetilde{\mathbf{G}}_{1}^{n-m})) \nonumber \\
&=&P(E_{n-m,k}) \ \textrm{ a.s.\ on%
}\ E_{m,k}.  \label{cond}
\end{eqnarray}

Let $r_{k}=[\frac{k}{p_{k}}]$ and define the stopping times $\tau _{i,k}$
recursively by $\tau _{0,k}=0$ and 
\begin{eqnarray*}
\tau _{i,k}(\omega )=\min \{n\in \mathbb{N}:n>\tau _{i-1,k}(\omega )\ \textrm{and}\
\omega \in E_{n,k}\}
\end{eqnarray*}
and $\tau _{i,k}(\omega )=\infty $ if $\omega $ lies in less than $i$ of the
sets $E_{n,k}.$ Lemma \ref{recurrevents} allows us to conclude from (\ref{cond}%
) and (\ref{expcount}) that with probability at least $1-\frac{1}{k}$ at
least $r_{k}$ of the sets $E_{n,k}$ occur, i.e. 
\begin{equation}
P(\tau _{r_{k},k}<\infty )\geq 1-\frac{1}{k}.  \label{highpr}
\end{equation}
We consider the events 
\begin{eqnarray*}
H_{i,k}=\{\tau _{i,k}<\infty \}\cap \{c_{k}\leq |X_{\tau _{i,k}+1}|,\textrm{%
\textrm{sgn}}(X_{\tau _{i,k}+1}) \neq \textrm{\textrm{sgn}}(S_{\tau _{i,k}}),%
\tilde{X}_{\tau _{i,k}+1}=0\}
\end{eqnarray*}
and $H_{k}=\bigcup_{i=1}^{r_{k}}H_{i,k}.$ 
Let $W_{n}=V_{n}\cup \{\omega :S_{n}S_{n+1}=0,\widetilde{S}_{n}=0=\widetilde{%
S}_{n+1}\}$, where $V_{n}$ is as in part (b) of Theorem \ref{Example}.
By definition $H_{i,k}\subset
E_{\tau _{i,k},k}$ and hence $H_{i,k}$ is contained in the set $W_{n}$ with $n=\tau _{i,k}$. So, if $\omega \in H_{k}$ then
there is some $n$ with $\omega \in W_{n}$ and $|X_{n+1}|\geq c_{k}$. Since
the sequence $(c_{k})$ is unbounded every point in $\limsup H_{k}$ lies in
infinitely many sets $W_{n}$. We show that 
\begin{equation}
\lim_{k\rightarrow \infty }P(H_{k})=1.  \label{gk}
\end{equation}
Denote by $U^{(k)}$ a random variable with law $R[0,y_{k}]$. By definitions
(\ref{Tdef}) and (\ref{Xdef}) we have 
\begin{eqnarray*}
P(c_{k}\leq |X|) &\geq &P(G=1)p_{k}P(c_{k}\leq U^{(k)}) \\
&\geq &\frac{2}{3}p_{k}(1-\frac{c_{k}}{y_{k}})>\frac{1}{2}p_{k}
\end{eqnarray*}
for all sufficiently large $k$ since $\frac{c_{k}}{y_{k}}\rightarrow 0$ because of (%
\ref{pest}). Let $\delta =P(X=0)$. Then the complements $H_{i,k}^{c}$ of our
sets satisfy 
\begin{eqnarray*}
P(\{\tau _{i,k}<\infty \}\cap H_{i,k}^{c}\ |\ \mathcal{F}_{\tau _{i,k}})&\leq&
1-P(c_{k}\leq |X|,sgn(X)>0,\tilde{X}=0) \\
&\leq& 1-\frac{\delta }{4}p_{k}
\end{eqnarray*}
and an induction argument shows that 
\begin{eqnarray*}
P\left( \{\tau _{r_{k},k}<\infty \}\cap
\bigcap_{j=1}^{r_{k}}H_{j,k}^{c}\right) \leq (1-p_{k}\frac{\delta }{4}%
)^{r_{k}}.
\end{eqnarray*}
The right-hand side in this inequality becomes arbitrarily small for large $k
$ by the choice of $r_{k}$. Hence by (\ref{highpr}) 
\begin{eqnarray*}
\lim_{k\rightarrow \infty }P(H_{k})=\lim_{k\rightarrow \infty }P(\tau
_{r_{k},k}<\infty )=1.
\end{eqnarray*}
We have shown that almost surely the event $W_{n}$ occurs for infinitely many $n$. By the transience of our random walk, the event $\{\omega :S_{n}S_{n+1}=0,\widetilde{S}_{n}=0=\widetilde{%
S}_{n+1}\}$ a.s. cannot occur infinitely often. So we conclude that a.s. for infinitely many $n$ the event $V_{n}$ occurs.
This completes the proof.

\section{Appendix: Some tools}

\begin{Lemma}[A counting variable estimate]
\label{recurrevents} Let $(\mathcal{F}_{n})_{0\leq n\leq N}$ be a (finite or
infinite) filtration. Let $(E_{n})$ be an adapted sequence of events such
that $E_{0}=\Omega $ and for $m<n$ 
\begin{eqnarray*}
P(E_{n}{\ |\ }\mathcal{F}_{m})\leq P(E_{n-m})\ a.s.\ on\ E_{m}.
\end{eqnarray*}
Let $\Phi $ be the number of events which occur (including $E_{0}$). Then
for each $r=0,1,2,\cdots $ 
\begin{equation}
P(\Phi >r)\geq 1-\frac{r}{\mathbb{E}(\Phi )}.
\end{equation}
\end{Lemma}

\begin{Proof}
Clearly it suffices to consider the case of finite $N$. We call an index $n$
a success time if $E_{n}$ occurs. Let $\tau _{r}(\omega )$ be the $r$-th
success time $\geq 1$ and let $\tau _{r}=N+1$ if $\Phi \leq r$. Moreover let 
$\Phi _{m}$ be the number of success times $\geq m$. Then the inequality in
our assumption implies for each $m\leq N$ 
\begin{eqnarray*}
\mathbb{E}(\Phi _{m}|\mathcal{F}_{m})=\sum_{n=m}^{N}P(E_{n}|\mathcal{F}%
_{m})\leq \sum_{n=m}^{N}P(E_{n-m})\leq \mathbb{E}(\Phi )
\end{eqnarray*}
(a.s. on $E_{m}$) and hence also $\mathbb{E}(\Phi _{\tau _{r}}|\mathcal{F}_{\tau _{r}})\leq 
\mathbb{E}(\Phi )$ on the set $\{\Phi >r\}=\{\tau _{r}\leq N\}\in \mathcal{F}%
_{\tau _{r}}$. Then 
\begin{eqnarray*}
\mathbb{E}(\Phi )\leq rP(\Phi \leq r)+\mathbb{E}\left( \mathbf{1}_{\{\Phi
>r\}}\left( \Phi _{\tau _{r}}+r\right) \right) \leq r+P(\Phi >r)\mathbb{E}%
(\Phi ).
\end{eqnarray*}
Dividing by $\mathbb{E}(\Phi )$ yields the result.
\end{Proof}

\begin{Lemma}
\label{momest} Let $T$ be a random variable on $\{0,1,2,...\}$ with
non-increasing weights. Then the estimate 
\begin{eqnarray*}
(\mathbb{E}T)^{2}\leq \frac{3}{4}\mathbb{E}T^{2}
\end{eqnarray*}
is valid. It implies $(\mathbb{E}T)^{2}\leq 3\mathrm{Var}(T).$
\end{Lemma}

\begin{Proof}
It is easy to see that a random variable on $\{0,1,2,...\}$ has
non-increasing weights iff it can be represented as a mixture of uniform
distributions $R[0,y]$ as in (\ref{Tdef}). So we get by Jensen's
inequality 
\begin{eqnarray*}
(\mathbb{E}T)^{2} &=&(\sum_{l=1}^{\infty }p_{l}y_{l}/2)^{2} \\
&\leq &\sum_{l=1}^{\infty }p_{l}(y_{l}/2)^{2}=\frac{3}{4}\sum_{l=1}^{\infty
}p_{l}y_{l}^{2}/3\leq \frac{3}{4}\sum_{l=1}^{\infty }p_{l}y_{l}(2y_{l}+1)/6=%
\frac{3}{4}\mathbb{E}T^{2}.
\end{eqnarray*}
\end{Proof}

\begin{Lemma}
\label{unimod}Let $T_{i},i\in \mathbb{N}$ be an identically distributed
sequence of random variables on $\{0,1,2,...\}$ with non-increasing weights.
Let $\epsilon \in {\pm 1}$ be a cointossing random variable and let $G$ be a
random variable with values in $\mathbb{N}$. If all these r.v.'s are
independent of each other then the law of 
\begin{equation}
X=\epsilon \sum_{i=1}^{G}(-1)^{i}T_{i}  \label{Xdef}
\end{equation}
is symmetric unimodal with 
\begin{equation}
\mathrm{Var}(T_{1})\leq \mathrm{Var}(X)\leq 4\mathbb{E}(G)\mathrm{Var}%
(T_{1}),  \label{varest}
\end{equation}
where the last estimate is interpreted trivially if $T_{1}$ has infinite
variance.
\end{Lemma}

\begin{Proof}
Denote by $\tau $ the law of the $T_{i}$ and by $\mu _{k}$ the law of $X$ in
the case where $G$ takes the constant value $k$. The convolution of $\tau $
with its reflected image on $-\mathbb{N}$ is symmetric and easily seen to be
unimodal. Moreover it is known that the convolution of two symmetric
unimodal laws is again symmetric unimodal. (Decompose both laws as mixtures
of uniform distributions on suitable centered intervals.) This implies the
assertion if $G$ is an even constant, i.e. $\mu _{2m}$ is symmetric
unimodal. Now assume that $G=2m+1$ is odd. The conditional law of $\epsilon
\sum_{i=1}^{m}T_{2i}-T_{2i-1}$ given $\epsilon $ is by symmetry equal to $%
\mu _{2m}$ and hence independent of $\epsilon $. Thus $\mu _{2m+1}$ is the
convolution of $\mu _{2m}$ with the law of $\epsilon T_{2m+1}$. The latter
is also symmetric unimodal and hence $\mu _{2m+1}$ is unimodal symmetric as
well. Since this property is also stable under mixtures the result follows
also for nonconstant $G$. Set $Y=\epsilon ^{-1}X$. We have
\begin{eqnarray*}
\mathrm{Var}(\epsilon Y) &=&\mathbb{E}(Y^{2})=\mathrm{Var}(Y)+(\mathbb{E}%
Y)^{2} \\
&=&\mathbb{E}G\cdot \mathrm{Var}(T_{1})+P(G\textrm{ odd})(\mathbb{E}%
T_{1})^{2}\leq 4\mathbb{E}(G)\mathrm{Var}(T_{1})
\end{eqnarray*}
where the variance of the sum $Y=\epsilon ^{-1}X$ is computed by Wald's
identity and Lemma \ref{momest} was used. This relation implies both the lower and upper estimate of the variance.
\end{Proof}

We follow the convention to call a
symmetric random $\mathbb{R}$-valued variable $T$ unimodal if it is absolutely continuous with respect to Lebesgue measure and the density can be chosen non-increasing on $\mathbb{R}_{+}$. Analogously, if $T$ is a symmetric random variable with values on $\mathbb{Z}$ it is called unimodal if it has non-increasing weights on $%
\{0,1,2,...\}$.

\begin{Lemma}
\label{unimodest} There is a positive constant $d$ such that for every
symmetric unimodal distribution $\mu $ with finite variance $\sigma ^{2} >0$
and every $c>0$ with $c\leq \sigma$ one has 
\begin{equation}
\mu (\{x:|x|<c\})\geq \frac{cd}{\sigma }.  \label{lest}
\end{equation}
\end{Lemma}

\begin{Proof}
1. First consider the case where $\mu $ is carried by $\mathbb{R}$.
Introduce a scaling parameter $\lambda \geq 1$ and observe that the
assumption on $\mu $ implies that for any $c>0$ we have $\mu (\{x:|x|<c\})\geq
\lambda ^{-1}\mu (\{x:|x|<\lambda c\})$. (Substitute $x$ by $\lambda x$ in the integral over the density function of $\mu$.) Hence we have with $\lambda
^{\prime }:=\lambda c/\sigma \geq c/\sigma $ chosen arbitrarily 
\begin{eqnarray*}
\mu (\{x :|x|<c\})\geq \lambda ^{-1}(1-\mu (\{x:|x|\geq \lambda c\})) \\
\geq \lambda ^{-1}(1-\frac{\sigma ^{2}}{\lambda ^{2}c^{2}})=\frac{c}{%
\sigma }(\lambda ^{\prime })^{-1}(1-(\lambda ^{\prime })^{-2})
\end{eqnarray*}
by Chebychev's inequality. Now we choose $d:=\max_{\lambda ^{\prime} > 1}((\lambda ^{\prime })^{-1}(1-(\lambda ^{\prime })^{-2}))$. This proves the assertion in the case that $\mu $
is carried by $\mathbb{R}$. 

2. Now consider the case that $\mu $ is carried by $\mathbb{Z}$. For $\sigma ^{2} < 3/4$ the it follows from Chebychev's inequality that $d$ can be chosen to be $1/4$. So we may assume $3/4 \leq \sigma ^{2}$ and even $3/4 \leq c^{2} \leq \sigma ^{2}$. Let $\nu $ be the uniform
distribution on $[-\frac{1}{2},\frac{1}{2}]$. It is easily checked that $%
\widehat{\mu }:=\mu \ast \nu $ is symmetric unimodal on $\mathbb{R}$ with
variance $\sigma ^{2}+\frac{1}{12}$ and $\mu (\{j\})=\widehat{\mu }([j-\frac{%
1}{2},j+\frac{1}{2}])$. We have for $3/4 \leq c^{2} \leq \sigma ^{2}$%
\begin{eqnarray*}
\mu (\{x :|x|<c\})\geq \widehat{\mu }(\{x:|x|<c-\frac{1}{2}\}) \\
\geq \frac{(c-\frac{1}{2})d}{\sqrt{\sigma ^{2}+\frac{1}{12}}}\geq \frac{%
cd^{\prime }}{\sigma }
\end{eqnarray*}
for a suitable $d^{\prime }>0$.

\end{Proof}

\begin{Lemma}
Let $\{\mu _{i}\}_{i=1}^{m}$ be a set of equidistributions on some integer
intervals of equal length $y \geq 1$, i.e. $\mu _{i}=R[a_{i},a_{i}+y]$.
Then we have
\begin{equation}
\max_{x\in \mathbb{Z}}(\mu _{1}\ast \mu _{2}\ast ...\ast \mu
_{m})(\{x\})\leq \frac{D}{\sqrt{m}y},  \label{maxest}
\end{equation}
where $D$ is some absolute constant.
\end{Lemma}

\begin{Proof} The result can be obtained by standard estimates for concentration
functions involving characteristic functions. The expression on the
left-hand side of (\ref{maxest}) is the concentration function $Q(\mu _{1}\ast
\mu _{2}\ast ...\ast \mu _{m};1/2)$ of the $m$-fold convolution, evaluated
for an interval length $\lambda =1/2$. We make use of an estimate given in
\cite{Salikhov}, which is essentially due to Esseen \cite{Esseen}:
\begin{equation}
Q(\mu _{1}\ast \mu _{2}\ast ...\ast \mu _{m};\lambda )\leq \lambda \left(
2\tau \left( \frac{\sin (\tau /2)}{\tau /2}\right) ^{2}\right)
^{-1}\int_{-2\tau /\lambda }^{2\tau /\lambda }|\varphi _{\mu _{1}\ast
...\ast \mu _{m}}(t)|dt  \label{charest}
\end{equation}
being valid for arbitrary $\lambda >0$ and $0<\tau <2\pi $. In order to
simplify the considerations, we substitute each $\mu _{i}$ by the same $%
\widehat{\mu }:=\mu _{i}((\cdot )+a_{i}+y/2)$ being symmetric with respect
to the origin. The resulting shift does not change the concentration
function of the convolution, but the characteristic functions become
real-valued. The explicit expression for the characteristic function $%
\varphi _{\widehat{\mu }}$ is given by $\frac{(y+1)^{-1}\sin \frac{1}{2}%
t(y+1)}{\sin \frac{1}{2}t},t\notin 2\pi \mathbb{Z}$. Choosing $\tau =\frac{1%
}{4(y+1)}$ we obtain, for some absolute constant $C$, from (\ref{charest})
\begin{eqnarray*}
Q(\mu _{1}\ast \mu _{2}\ast ...\ast \mu _{m};\frac{1}{2}) &\leq &C\frac{1}{y}%
\int_{0}^{(y+1)^{-1}}\left| \frac{\sin \frac{1}{2}t(y+1)}{(y+1)\sin \frac{1}{%
2}t}\right| ^{m}dt.
\end{eqnarray*}
Now we make use of the fact that the members in the Taylor expansion of the sine
function are alternating and non-increasing in the interval considered. We
may continue the estimate as follows
\begin{eqnarray*}
&\leq &C\frac{1}{y}\int_{0}^{1}\left( \frac{\sin \frac{1}{2}t}{(y+1)\sin \frac{1}{2}t(y+1)^{-1}}\right) ^{m}dt \\
&\leq &C\frac{1}{y}\int_{0}^{1}\left( \frac{\frac{1}{2}t-\frac{1}{12}t^{3}+%
\frac{1}{240}t^{5}}{\frac{1}{2}t-\frac{1}{12(y+1)^{2}}t^{3}}\right) ^{m}dt \\
&\leq &C\frac{1}{y}%
\int_{0}^{1}\left( \frac{1-\frac{1}{6}t^{2}+\frac{1}{120}t^{4}}{1-\frac{1}{%
6(y+1)^{2}}t^{2}}\right) ^{m}dt \\
&\leq &C\frac{1}{y}\int_{0}^{1}\left( \frac{1-\frac{1}{5}t^{2}}{1-\frac{1}{%
6(y+1)^{2}}t^{2}}\right) ^{m}dt \\
&\leq &C\frac{1}{y}\int_{0}^{1}\left( \frac{1-\frac{1}{5}t^{2}}{1-\frac{1}{10%
}t^{2}}\right) ^{m}dt \\
&\leq &C\frac{1}{y}\int_{0}^{1}\left( 1-\frac{1}{10}t^{2}\right) ^{m}dt \\
&\leq &C\frac{1}{y}\int_{0}^{1}e^{-\frac{1}{10}t^{2}m}dt \\
&\leq &C\frac{1}{y}\int_{0}^{\infty }e^{-\frac{1}{10}t^{2}m}dt \\
&=&\frac{D}{\sqrt{m}y},
\end{eqnarray*}
for some absolute constant $D$.
\end{Proof}

\textbf{Acknowledgement}. The authors are very grateful to the referee 
for his many highly constructive and painstaking comments to the first 
version of this paper. In particular they helped with some pitfalls in the 
somewhat involved construction above.

\end{document}